\documentclass[11pt]{article}
\usepackage{amssymb,amsfonts,amsthm}
\usepackage{color}
\usepackage[latin1]{inputenc}
\usepackage[active]{srcltx}
\def\openC{{\rm C\kern-.18cm\vrule width.8pt height 7pt depth-.2pt \kern.18cm}}
\def\openN{{{\rm I}\kern-.16em {\rm N}}}
\def\openR{{{\rm I}\kern-.16em {\rm R}}}
\def\openT{{{\rm T}\kern-.42em {\rm T}}}
\def\openZ{{{\rm Z}\kern-.28em{\rm Z}}}

\def\eop{\hfill\rule{2.5mm}{2.5mm}}
\def\pf{\par\smallbreak\noindent {\bf Proof.} \ }

\newtheorem{thm}{Theorem}[section]

\newtheorem{lem}[thm]{Lemma}
\newtheorem{prop}[thm]{Proposition}

\theoremstyle{definition}

\newtheorem{ex}[thm]{Example}

\def\eop{\hfill\rule{2.5mm}{2.5mm}}
\begin{document}

\title{
{\textbf{\Large{Jackson kernels: a tool for analyzing the decay of eigenvalue sequences of integral operators on the sphere}}} \vspace{-4pt}
\author{
T. Jord\~{a}o\,\thanks{
Both authors partially supported by FAPESP, grants $\#$ 2014/06209-1 and $\#$ 2014/00277-5 respectively. }
\,\,\&\,\,
V. A. Menegatto}
}

\date{}
\maketitle \vspace{-30pt}
\bigskip

\begin{center}
\parbox{13 cm}{{\small Decay rates for the sequence of eigenvalues of positive and compact integral operators has been largely investigated for a long time in
the literature.\ In this paper, the focus will be on positive integral operators acting on square integrable functions on the unit sphere and generated by a
kernel satisfying a H\"{o}lder type assumption defined via average operators.\ In the approach to be presented here, the decay rates will be reached from convenient
estimations on the eigenvalues of the operator themselves, with the help of specific properties of a generic approximation operator defined through the so-called
generalized Jackson kernels.\ The decay rates have the same structure of those known to hold in the cases in which the H\"{o}lder condition is the classical one.\
Therefore, within the spherical setting, the abstract approach to be introduced here extends some classical results on the topic.
}}
\end{center}



\thispagestyle{empty}

%
%

\section{Introduction}\label{s1}

The present paper provides a concise approach to obtain decay rates for the eigenvalue sequence of positive integral operators acting on square
integrable functions on the sphere, in the case when the generating kernel of the operator satisfies an abstract H\"{o}lder condition.\ A brief feedback about this subject
cannot omit results obtained in the late 80's by several authors.\ This section begins with a description of some of those results, mainly those which pertain to the scope of this paper.

We start with a function $K$ in $L^2([0,1]^2)$ and consider the compact operator ${\cal{L}}_K: L^2([0,1])\rightarrow L^2([0,1])$ generated by it
$$
{\cal{L}}_K(f)(x)=\int_{0}^1 K(x,y)f(y)\,dy,  \quad f \in L^2([0,1]),\quad x \in [0,1].
$$
In this case, and also in others in which the interval is replaced with a more general space, we simply call $K$ the kernel and ${\cal{L}}_K$ the operator.\ The introduction of the symmetry assumption
$$
K(x,y) = \overline{K(y,x)}, \quad (x,y) \in [0,1]^2,
$$
makes the operator ${\cal{L}}_K$ self-adjoint and,
therefore, its eigenvalue sequence $\{\lambda_n\}$ can be ordered in a decreasing manner
$$
|\lambda_1| \geq |\lambda_2| \geq |\lambda_3| \geq \cdots,.
$$
taking into account multiplicities.\ In particular, $\{\lambda_n\}$ approaches zero as $n \to \infty$.\ By the way, the equality above and some others in the paper need to be interpreted as equalities a. e..

A classical result of Weyl states that
$$
\lambda_n = o(n^{-k-1/2}), \quad (n\to \infty),
$$
whenever $K \in C^k([0,1]^2)$.\ After the introduction of a positivity assumption, Reade (\cite{reade-1}) established the faster decay rate
$$
\lambda_n = o(n^{-k - 1}), \quad (n\to \infty).
$$
The positivity mentioned above refers to the property
$$
\int_0^1 \int_0^1 K(x,y)f(x) \overline{f(y)}\,d x dy \ge 0, \quad f \in L^2([0,1]),
$$
which, in the most important cases, corresponds to the usual positive definiteness of the kernel $K$.

Later on, still keeping the positiveness in the setting, the same author deduced the decay rate (\cite{reade-2})
$$
\lambda_n = O(n^{-1-r}), \quad (n \to \infty).
$$
under the following H\"{o}lder assumption on the generating kernel $K$:
$$
|K(x,y) - K(z,w)| \le C(|x-z|^r +|y-w|^r ), \quad x,y,z,w \in [0,1],\quad  r\in (0,1).$$

A few years later, an outstanding generalization of the results above appeared in \cite{kuhn2} with the replacement of $[0,1]$
with a compact $C^{\infty}$ manifold.\ If the manifold is the usual $m$-dimensional unit sphere $S^m$ endowed with its surface measure
$\sigma_m$, the compact operator ${\cal{L}}_K$ now acts on $L^2(S^m):=L^2(S^m,\sigma_m)$, that is,
$$
{\cal{L}}_K(f)(x)=\int_{S^m} K(x,y)f(y)\,d\sigma_m(y),  \quad f \in L^2(S^m),\quad x \in S^m,
$$
in which $K \in L^2(S^m \times S^m):=L^2(S^m \times S^m,\sigma_m\times \sigma_m)$.\ If $K$ is continuous and satisfies the standard H\"{o}lder condition
$$
|K(x,w)-K(y,w)| \leq B(w)d_m(x,y)^{\beta},\quad x,y,w \in S^m,
$$
for some $B\in L^1(S^m)$ and $\beta\in [0,1)$, the main result in \cite{kuhn2} produced the decay
$$
\lambda_n = O(n^{-1-\beta/m}), \quad (n \to \infty),
$$
for the eigenvalue sequence $\{\lambda_n\}$ of ${\cal{L}}_K$.

The papers mentioned above included examples of integral operators for which the decay of the eigenvalue sequence matches exactly the decay obtained.\ In other
words, in all cases above, the decay rates are best possible within each setting.

The setting in the present paper will be spherical one and the focus will be on integral operators generated by (not necessarily continuous) kernels satisfying a H\"{o}lder type assumption
defined by average operators.\ Below we stress a few points the reader should consider before and during the reading of this paper.\ They provide a reason why we have written the paper
in the format it is:\\
$(i)$ The H\"{o}lder assumptions we will employ are weaker then those found in other references, even when the spherical setting is considered;\\
$(ii)$ The spherical setting allows many different H\"{o}lder conditions which, to the best of our knowledge, are not equivalent to each other;\\
$(iii)$ The spherical setting allows different approaches to the problem and that permits variations in the assumptions;\\
$(iv)$ The approach adopted here allows the establishment of an abstract setting under which the decay rates can searched;\\
$(v)$ The approach adopted here is comparable to others found in the literature, however, it has its own characteristics;\\
$(vi)$ The approach permits the consideration of slightly weaker general assumptions still reaching the same decay rates found in the literature; \\
$(vi)$ All the results to be proved can be considered in more settings, as long as it has a background structure similar to that available in the spherical setting (two-point homogeneous spaces for example);\\
$(viii)$ The spherical setting has practical relevancy in other areas, for instance, in Geo-mathematics and meteorological sciences in general.

An outline of the paper is as follows.\ In Section 2, we introduce the abstract setting along with an abstract H\"{o}lder condition defined by the spherical convolution operator on which the
main results of the paper will be based upon.\ We also include two motivational examples that may justify the why we consider a general and abstract setting.\
In Section 3, we define the so called approximation operators, here manufactured with the help of the generalized Jackson kernels.\ At the end of the section we show the approximation operator has finite rank when the setting is either one of the two motivational examples.\ Section 4 begins with the notion of positive integral operator that pertains to this work.\ It is followed by a technical estimate for integrals involving the generalized Jackson kernels and inequalities involved in the estimation of the approximations numbers of the square root of the positive integral operator.\ That is followed by inequalities for the approximation numbers of the integral operator itself, under the assumption that the rank of the attached approximation operator is finite.\ Finally, the section is closed with the main results in the paper.\ Section 5 is reserved for relevant remarks and the pointing of some open questions.


\section{A H\"{o}lder condition based on spherical convolutions}

Let us begin with the basic structure to be used in the paper.\ In addition to the spaces $L^2(S^m)$, we will stick to the usual spaces $L^p(S^m)$, $p=1,\infty$.\ The norm in all of then will be written $\|\cdot \|_p$, $p=1,2,\infty$.\ Finally, we want to emphasize from the outset that, throughout the whole paper, that the dimension $m$ will be fixed.

We will consider a H\"{o}lder assumption based upon a fixed family of nonnegative functions $\{\mathcal{Z}_t^m: t\in (0,\pi)\}$ belonging to $L^1([-1,1], d\omega_m)$, in which
$$
d\omega_m(u)=(1-u^2)^{(m-2)/2}du,\quad u \in [-1,1].
$$
If $\tau_m$ is the surface measure of $S^m$, then the norm in this space is
$$
\|\phi\|_{1,m}:=\frac{\tau_{m-1}}{\tau_m}\int_{-1}^{1}|\phi(u)|d\omega_m(u),\quad \phi \in L^1([-1,1], d\omega_m),
$$
and the formula
$$
Z_t^m(x,y):=\mathcal{Z}_t^m(x \cdot y),\quad x,y \in S^m,
$$
defines an associated family $\{Z_t^m: t\in (0,\pi)\}$ of isotropic kernels on $S^m$.\ As usual, {\em isotropy} of a kernel refers to its invariance with respect to orthogonal transformations
of the space where $S^m$ sits.\ The setting to be undertaken here demands two assumptions:\\
${\bf A1}$ - The family $\{\mathcal{Z}_t^m: t\in(0,\pi)\}$ is uniformly bounded in $L^1([-1,1], d\omega_m)$.\\
${\bf A2}$ - If $V_m(t)$ is the surface area of the support of $Z_t^m(x, \cdot)$ in $S^m$, then there exist a positive integer $\alpha(m)$ and positive constants $c_{m}$ and $C_{m}$ so that
$$ V_m(t)\leq C_{m}t^{\alpha(m)},\quad t\in (0,\pi),$$
and
$$c_m t^{\alpha(m)}\leq V_m(t),\quad t\in (0,\pi/2).$$

The surface area of the support of $Z_t^m(x, \cdot)$ in $S^m$ mentioned in $A2$ does not depend upon $x \in S^m$.\ Indeed, fix $x_1,x_2\in S^m$ and for each $i\in \{1,2\}$, write $S_i^m(t)$ to denote the support of $Z_t^m(x_i,\cdot)$, and put
$$V_m(t,x_i):=\int_{S_i^m(t)}d\sigma_m(z).$$
Since each $Z_t^m$ is not necessarily continuous, we have that
$$S_i^m(t)=S^m -\cup \{A : A\ \mbox{is open in\,} S^m\ \mbox{and\ } \sigma_m(\{x \in A: Z_t^m(x_i,x)\neq 0\})=0\}, \quad i=1,2.$$
Using the isotropy of $Z_t^m$ and some straightforward computations, it is easily seen that $\mathcal{O}(S_1^m(t))=S_2^m(t)$, whenever $\mathcal{O}$ is an orthogonal transformation of $\mathbb{R}^{m+1}$ satisfying $\mathcal{O}(x_1)=x_2$.\ It is now clear that
$$V_m(t,x_2)=\int_{\mathcal{O}(S_{1}^m(t))}d\sigma_m(z)=\int_{S_{1}^m(t)}d\sigma_m(z)=V_m(t,x_1).$$

Under the setting introduced above, we define
$$
T_t(f):= Z_t^m *f, \quad f \in L^p(S^m), \quad t\in (0,\pi), \quad p=1,2,\infty
$$
in which
$$(Z_t^m * f)(x)=\frac{1}{\tau_m}\int_{S^m} \mathcal{Z}_t^m(x \cdot y)f(y)d\sigma_m(y), \quad x \in S^m, \quad f \in L^p(S^m),$$
is the {\em spherical convolution} of $Z_t^m$ and $f$ in $S^m$.\ Every $T_t$ is a well-defined bounded linear operator from $L^p(S^m)$ into itself with $\|T_t\| \leq \|\mathcal{Z}_t^m\|_{1,m}$.\
If $\beta \in (0,2]$ and $B$ is a nonnegative function from $L^{\infty}(S^m)$, then a kernel $K$ on $S^m$ is {\em $(T_t,B,\beta)$-H\"{o}lder} if
$$
|T_t(K(y,\cdot))(x)-K(y,x)|\leq B(y)t^{\beta}, \quad x,y\in S^m,\quad  t\in (0,\pi).
$$
Below, we discuss two particular cases which served as motivation for the abstract setting introduced above
and also for the consideration of the H\"{o}lder condition just defined.

\begin{ex} For $n\geq 1$ and $t \in [0,\pi]$, let
$C_m(t)$ be the total volume of the cap $$C_t^x=\{y\in S^m: x \cdot y \geq \cos t\}$$ of $S^m$ defined by
$t$ and ``the pole" $x$.\ Clearly,
$$
C_m(t)=\tau_{m-1}\int_0^{t}(\sin h)^{m-1}dh, \quad t \in (0,\pi),
$$
a quantity that does not depend upon $x$.\ The formula
$$\mathcal{Z}_{n,t}^{m}(x\cdot y)=\left\{\begin{array}{rc}
\tau_m C_m(t)^{-1}(x\cdot y-\cos t)^{n-1}(1-\cos t)^{-(n-1)},& \quad \mbox{if\ \ } \cos t \leq x \cdot y\leq 1\\
0,& \quad otherwise
\end{array}\right.$$
defines families $\{\mathcal{Z}_{n,t}^{m}\}$ of locally supported kernels on $S^m$.\ A construction of such kernels by iteration with spherical convolutions can be found in \cite{levesley,schreiner}.\ Since $$\|\mathcal{Z}^m_{n,t}\|_{1,m}=\frac{\tau_{m-1}}{\tau_m}\int_{-1}^1 \mathcal{Z}^m_{n,t}(u)d\omega_m(u)\leq\frac{\tau_{m-1}}{C_m(t)}\int_{\cos t}^1 (1-u^2)^{(m-2)/2}du=1,$$
each family $\{\mathcal{Z}^m_{n,t} : t \in (0,\pi)\}$ is uniformly bounded in $L^1([-1,1], d\omega_m)$.\ On the other hand, it is easily seen that the surface area $V_m(t)$ of the support of $Z_{n,t}^m(x, \cdot)$ in $S^m$ is precisely $C_m(t)$ (there is no dependence on $n$) while direct computation yields
$$\frac{\tau_{m-1}}{m}\left(\frac{2}{\pi}\right)^{m-1} t^m\leq V_m(t) \leq \tau_{m-1}t^m,  \quad t \in (0,\pi).$$
Thus, both assumptions A1 and A2 hold in this case.\ The particular case $n=1$ recovers the usual {\em average operator} $M_t$ on $S^m$ (\cite{berens}), that is,
$$
M_t(f)(x)=(\mathcal{Z}_{1,t}^m * f)(x)=\frac{1}{C_m(t)}\int_{C_t^x}f(w)dr(w), \quad x\in S^m,\quad  t\in (0,\pi),
$$
and the abstract $(T_t,B,\beta)$-H\"{o}lder condition turns itself into the {\em averaged H\"{o}lder condition}
$$
|M_t(K(y,\cdot))(x)-K(y,x)|\leq B(y)t^{\beta}, \quad x,y\in S^m,\quad  t\in (0,\pi).
$$
\end{ex}

\begin{ex} Here we will consider the \emph{Stekelov-type mean} operator introduced and discussed in \cite{ditzian2}.\ If $R_m(t):= \tau_{m-1}(\sin t)^{m-1}$, $t \in (0,\pi)$, then it has the form
$$
E_t(f)(x)=\frac{1}{D_m(t)}\int_0^{t}\frac{C_m(s)}{R_m(s)} M_s(f)(x)ds, \quad x\in S^m,\quad  t\in (0,\pi),
$$
the normalizing constant $D_m(t)$ being chosen so that $E_t(1)=1$.\ In order to see that the operators $E_t$ fit into the convolution structure we are using,
it suffices to consider the family of locally supported kernels
$$W_t^m(x,y):=\mathcal{W}_{t}^m(x \cdot y):=\left\{\begin{array}{rc}
\displaystyle{\int_0^t \frac{1}{R_m(s)}\mathcal{Z}_{1,s}^m(x \cdot y)ds},& \quad \mbox{if\ \ } \cos t \leq x\cdot y\leq  1\\
0,& \quad \mbox{otherwise},
\end{array}\right.$$
where $\mathcal{Z}_{1,s}^m$ are the kernels described in the previous example.\ Clearly,
\begin{eqnarray*}
\|\mathcal{W}^m_{t}\|_{1,m} & = & \frac{\tau_{m-1}}{\tau_m}\int_{-1}^1 \mathcal{W}^m_{t}(u)d\omega_m(u)\\
& \leq & \frac{\tau_{m-1}}{\tau_m D_m(t)}\int_{\cos t}^1 \left[\int_0^t \frac{C_m(s)}{R_m(s)}|\mathcal{Z}_{1,s}^m(u)|ds\right](1-u^2)^{(m-2)/2}du\\
& = & \frac{1}{D_m(t)}\int_0^t  \frac{C_m(s)}{R_m(s)}\left[\frac{\tau_{m-1}}{\tau_m}\int_{\cos t}^1 |\mathcal{Z}_{1,s}^m(u)(1-u^2)^{(m-2)/2}du\right]ds\\
& = & \frac{1}{D_m(t)}\int_0^t  \frac{C_m(s)}{R_m(s)}ds.
\end{eqnarray*}
The normalization we have chosen for the $D_m(t)$ provides the uniform boundedness for the family $\{\mathcal{W}^m_{t}: t \in (0,\pi)\}$.\ The support of $W_t^m(x,\cdot )$ is $C_m(t)$ and $A2$
holds as in the first example.\ Since the kernel $\mathcal{W}_{t}^m$ is isotropic, we have
$$
E_t(f)(x)=(W_{t}^m\ast f)(x), \quad f \in L^2(S^m),\quad t \in (0,\pi).
$$
In this case, the abstract $(T_t,B,\beta)$-H\"{o}lder condition turns itself into the {\em Stekelov-mean H\"{o}lder condition}
$$
|E_t(K(y,\cdot))(x)-K(y,x)|\leq B(y)t^{\beta}, \quad x,y\in S^m,\quad  t\in (0,\pi).
$$
\end{ex}

Before closing the section, let us return to the standard H\"{o}lder condition introduced in Section 1, now considering $B \in L^\infty(S^m)$.\ It is straightforward to
verify that a kernel $K$ satisfying the usual H\"{o}lder condition also satisfies an averaged H\"{o}lder condition (with the same index $\beta$ but not necessarily the same $B$).\ Likewise, a kernel satisfying an averaged H\"{o}lder condition also satisfies an Stekelov-mean H\"{o}lder condition (with the same index $\beta$ but not necessarily the same $B$).\ Thus, we have a chain of conditions from the stronger usual H\"{o}lder condition to the weaker Stekelov-mean H\"{o}lder condition.


\section{The approximating operators}

In this section, we introduce the approximation operators we intend to use in some critical arguments in the paper where we need to estimate the approximation numbers of our operators.\ They
will depend on the setting introduced in Section 2 and on the generalized Jackson kernels.\ The use of these kernels were influenced by the paper \cite{dikman} wherein standard Jackson kernels were used to
 obtain decay rates for the sequence of eigenvalues of the integral operator on $L^2([0,1]^2)$ in the case $K$ is differentiable in $[0,1]$ up to a certain order.\ On the other hand, it is well known that
 the generalized Jackson kernels imply optimal results in many problems in analysis and approximation theory.

We will assume that the setting at the beginning of Section 2 has been fixed here.\ For positive integers $l$ and $\mu\geq 2$, tied to each other via the formula $\nu=l(\mu-1)$, the {\em generalized Jackson kernel} is given by
$$
J_{\nu,n}(t):= \frac{1}{c_{\nu,n}}\left[\frac{\sin(\mu t/ 2)}{\sin(t/2)}\right]^{2l}, \quad t\in (0,\pi),\quad n\in \mathbb{Z}_+,
$$
with the normalization constant $c_{\nu,n}$ computed by the formula
$$c_{\nu,n}=\int_{0}^{\pi} \left[\frac{\sin(\mu t/ 2)}{\sin(t/2)}\right]^{2l} V_m(t) (\sin t)^n\,dt.$$
Here, the constant $V_m$ is that one introduced in A2.\ Clearly, the constants $c_{\nu,n}$ depend upon $m$ too, but that will be not enforced in the notation adopted.\ On the other hand, it is easily seen that the normalization corresponds to
$$
\int_{0}^{\pi}J_{\nu,n}(t)V_m(t) (\sin t)^n \, dt=1.
$$

The integral operators themselves can now be defined through the convolution operators $\{T_t:t\in(0,\pi)\}$.

\begin{prop} \label{phii} For each $n \in \mathbb{Z}_+$, the formula
$$\Phi_{\nu,n}(f)(x)=\int_0^\pi J_{\nu,n}(t)T_t(f)(x)V_m(t)(\sin t)^n dt, \quad f \in L^2(S^m), \quad x \in S^m,$$
defines a bounded linear operator $\Phi_{\nu,n}$ from $L^2(S^m)$ into itself.
\end{prop}
\pf Minkowski's inequality for integrals (\cite[p.194]{folland}) implies that
$$
\|\Phi_{\nu,n}(g)\|_2  \leq  \int_0^{\pi}J_{\nu,n}(t)\|T_{t}(g)\|_2V_m(t)(\sin t)^n\, dt, \quad g \in L^2(S^m).$$
Since $T_t$ is a convolution operator, it follows that
$$
\|\Phi_{\nu,n}(g)\|_2  \leq   \|g\|_2 \int_0^{\pi}\|\mathcal{Z}_t^m\|_{1,m}J_{\nu,n}(t)V_m(t)(\sin t)^n\, dt\leq M \|g\|_2 , \quad g \in L^2(S^m),$$
in which $M$ is a uniform bound for the sequence $\{\mathcal{Z}_t^m: t\in(0,\pi)\}$.\eop

In many cases, the formula
$$
a_n({\cal{L}}_K)=\min\{\|{\cal{L}}_K-U\|: \rho(U)\leq n-1\},
$$
in which $\rho(U)$ is the rank of $U$, is a useful tool in either the exact computation or the estimation of the $n$-th approximation number $a_n({\cal{L}}_K)$ of the operator of ${\cal{L}}_K$.\ As a matter of fact, $a_n({\cal{L}}_K)$ coincides with the $n$-th eigenvalue of the operator in those situations.\ So, if $\rho(\Phi_{\nu,n}) <\infty$, it is clear that the composition $U=\Phi_{\nu,n} \circ {\cal{L}}_K$ is eligible to be used in the estimation of some of approximation numbers.\ Since $J_{\nu,n}$ is an even trigonometric polynomial of degree $\nu$ (\cite{lizor}), it is reasonable to expect that $\rho(\Phi_{\nu,n}) <\infty$ for some special choices of $T_t$.\ The results that close this section will ratify that in the examples presented in the second half of Section 2.

We will write $\mathcal{H}_k^m$ to denote the space of all spherical harmonics of degree $k$ in $m+1$ variables and will denote its dimension by $N(m,k)$.\ The orthogonal decomposition $L^2(S^m)=\oplus_{k=0}^{\infty}\mathcal{H}_k^m$ is well-known while the orthogonal projection of $L^2(S^m)$ over $\mathcal{H}_k^m$ is given by the formula
$$
\mathcal{Y}_{k}(g)(x)=\frac{N(m,k)}{\tau_m P_k^{(m-1)/2}(1)}\int_{S^m} P_k^{(m-1)/2}(x\cdot y)g(y)d\sigma_m(y), \quad g \in L^2(S^m),\quad x\in S^m,
$$
in which $P_k^{(m-1)/2}$ is the usual Gegenbauer polynomial of degree $k$ associated to the dimension $m$.\
The additional formula
$$\mathcal{Y}_k(M_{t}(g))=\frac{\tau_{m-1}}{C_m(t)P_k^{(m-1)/2}(1)}\left(\int_{0}^{t}P_k^{(m-1)/2}(\cos h)(\sin h)^{m-1}dh\right)\mathcal{Y}_k(g), \quad g \in L^2(S^m), $$
for the projections of the elements in the range of $M_{t}$ was derived in \cite{berens}.\ A nice reference for the results on the analysis on the sphere mentioned above and ahead is \cite{daixu}.

The propositions below provide estimates for the rank of the operator in Proposition \ref{phii}, in the cases in which $T_t$ is either the average operator $M_t$ or the Stekelov-type mean operator $E_t$.

\begin{prop}\label{operatorthm} The operator
$$
\Phi_{\nu,1}(f)(x)=\int_0^\pi J_{\nu,1}(t)M_t(f)(x)C_m(t)\sin t \,dt, \quad f \in L^2(S^m), \quad x \in S^m,
$$
has rank at most $N(m+1,\nu+1)$.
\end{prop}

\pf For $f \in L^2(S^m)$ fixed, we will show that $\Phi_{\nu,1}(f)$ is a spherical
polynomial of degree at most $\nu+1$.\ That will imply the estimate announced in the statement of the proposition due to two facts: the space of all spherical polynomials of degree at most $\nu+1$ is precisely
$\oplus_{k=0}^{\nu+1}\mathcal{H}_k^m$ and its dimension is $\sum_{k=0}^{\nu+1}N(m,k)=N(m+1,\nu+1)$.\
The proof will be complete as long as we show that $\mathcal{Y}_{k}(\Phi_{\nu,1}(f))=0$, $k=\nu+2, \nu+3, \ldots$.\ Direct computation reveals that
$$
\mathcal{Y}_k(\Phi_{\nu,1}(f))(x)=\int_{0}^{\pi}J_{\nu,1}(t)\mathcal{Y}_k(M_{t}(f))(x)C_m(t)\sin t\, dt, \quad x \in S^m,
$$
while the formula prior to the statement of the theorem leads to
$$\label{projoperator2}
\mathcal{Y}_k(\Phi_{\nu,1}(f))=\frac{\tau_{m-1}}{P_k^{(m-1)/2}(1)}\left\{\int_{0}^{\pi}J_{\nu,1}(t)\left[\int_{0}^{t}P_k^{(m-1)/2}(\cos h)(\sin h)^{m-1}dh\right]\sin t\, dt\right\}\mathcal{Y}_k(f).
$$
The inner integral can be put into the form
$$
\int_{0}^{t}P_k^{(m-1)/2}(\cos h)(\sin h)^{m-1}dh=-\int_{1}^{\cos t}P_k^{(m-1)/2}(u)(1-u^2)^{(m-2)/2}du.
$$
Invoking the classical equality (\cite[p.81-82]{szeg})
$$
\frac{d}{du}\left[-\frac{m-1}{k(k+m-1)}(1-u^2)^{m/2}P_{k-1}^{(m+1)/2}(u)\right]= (1-u^2)^{(m-2)/2}P_{k}^{(m-1)/2}(u),
$$
we deduce that
$$
\int_{0}^{t}P_k^{(m-1)/2}(\cos h)(\sin h)^{m-1}dh=\frac{m-1}{k(k+m-1)}(\sin t)^{m}P_{k-1}^{(m+1)/2}(\cos t).
$$
Consequently,
$$
\mathcal{Y}_k(\Phi_{\nu,1}(f))=\frac{\tau_{m-1}(m-1)}{P_k^{(m-1)/2}(1)k(k+m-1)}\left[\int_{0}^{\pi}J_{\nu,1}(t)P_{k-1}^{(m+1)/2}(\cos t)(\sin t)^{m+1} dt\right]\mathcal{Y}_k(f).
$$
Since $J_{\nu,1}(t)$ is a polynomial of degree $\nu$ with respect to $\cos t$, it is easily seen that we can write it in the form
$$
J_{\nu,1}(\cos t)=\sum_{j=0}^{\nu}a_jP_j^{(m+1)/2}(\cos t),\quad a_1, a_2, \ldots, a_{\nu} \in \mathbb{R}.
$$
In particular,
$$
\int_{0}^{\pi}J_{\nu,1}(t)P_{k-1}^{(m+1)/2}(\cos t)(\sin t)^{m+1} dt=\sum_{j=0}^{\nu}a_j\int_{0}^{\pi}P_j^{(m+1)/2}(\cos t)P_{k-1}^{(m+1)/2}(\cos t)(\sin t)^{m+1} dt,
$$
and the well-known orthogonality relation (\cite[p.98]{morimoto})
$$\int_0^{\pi}P_k^{(m+1)/2}(\cos t)P_{l}^{(m+1)/2}(\cos t)(\sin t)^{m+1}dt=0, \quad k\neq l,
$$
implies that
$$\int_{0}^{\pi}J_{\nu,1}(t)P_{k-1}^{(m+1)/2}(\cos t)(\sin t)^{m+1} dt=0, \quad k-1\geq \nu+1.$$
It is now clear that
$$
\mathcal{Y}_k(\Phi_{\nu,1}(f))=0, \quad k= \nu+2, \nu+3, \ldots,
$$
and the result follows.\eop \vspace*{2mm}

\begin{prop}\label{operatorthm1} The operator
$$
\Phi_{\nu,m}(f)(x)=\int_0^\pi J_{\nu,n}(t)E_t(f)(x)D_m(t)(\sin t)^m dt, \quad f \in L^2(S^m), \quad x \in S^m,
$$
has rank at most $N(m+1,\nu+1)$.
\end{prop}
\pf Since it is analogous to the proof of the previous proposition, the details will be not included.\eop


\section{Decay rates via the generalized Jackson kernels}

The attention in this section will be directed to
integral operators ${\cal{L}}_K$ of the form
$$
{\cal{L}}_K(f)(x)=\int_{S^{m}}K(x,y)f(y)\,d\sigma_m(y), \quad x \in S^{m}, \quad f \in L^2(S^m).
$$
that possess the features below:\\
- it is generated by an element $K$ of $L^2(S^m\times S^m)$ (so, it is a linear operator from $L^2(S^m)$ into itself);\\
- the kernel $K$ is $L^2(S^m)$-positive definite in the sense that
$$
\int_{S^m}{\cal{L}}_K(f)(x)\overline{f(x)}d\sigma_m(x)\rangle \geq 0, \quad f \in L^2(S^m);
$$
- the square root ${\cal{L}}_K^{1/2}$ of ${\cal{L}}_K$ is an integral operator on $L^2(S^m)$ generated by a hermitian kernel $K_{1/2}:S^m \times S^m \to \mathbb{C}$;\\
- the original kernel $K$ can be recovered from $K_{1/2}$, that is,
\begin{equation}\label{recover}
\int_{S^m} K_{1/2}(x,y)K_{1/2}(w,x) d\sigma_m(x) = K(w,y), \quad y,w \in S^m.
\end{equation}
A usual concrete setting in which all the conditions above hold is described in \cite{schaback}.

An operator as above will be called a {\em positive integral operator} from now on.\ The category of positive integral operators includes those integral operators generated
by a continuous and positive definite kernel in the usual sense, as one can ratify in \cite{ferr}.\ A positive integral operator has countably many nonnegative eigenvalues which can be ordered as
$$
\lambda_1({\cal{L}}_K) \geq \lambda_2({\cal{L}}_K) \geq \cdots \geq 0,
$$
repetitions being included in accordance with algebraic multiplicities.\ After we order the eigenvalues of ${\cal{L}}_K^{1/2}$ in the same way, it holds
\begin{equation}\label{squarelambda}
\lambda_n({\cal{L}}_K^{1/2}) = (\lambda_n({\cal{L}}_K))^{1/2}=a_n({\cal{L}}_K))^{1/2}, \quad n=1,2,\ldots.
\end{equation}
and
$$
\|{\cal{L}}_K\| \geq a_1({\cal{L}}_K) \geq a_2({\cal{L}}_K)\geq \cdots \geq 0.
$$
For a general treatment on approximation numbers of operators, we refer the reader to \cite{piet} while a treatment in a setting similar to the one used here can be found in \cite{gohb}.

\subsection{The essential estimates}

This subsection contains preliminary estimates for the norm of the operator
$$\mathcal{L}_K^{1/2}-\Phi_{\nu}(\mathcal{L}_K^{1/2}): L^2(S^m)\rightarrow L^{\infty}(S^m)$$
when $\mathcal{L}_K$ is a positive integral operator generated by a $(T_t,B,\beta)$-H\"{o}lder kernel $K$.\ In particular, we remind the reader that the setting described in Section 2 needs to hold here,
including the assumptions A1 and A2.

Two reasons justify why we will estimate $a_n({\cal{L}}_K^{1/2})$ instead of $a_n({\cal{L}}_K)$: formula (\ref{squarelambda}) is available in the most important cases and the applications on decay rates for eigenvalues we seek demand the approximation numbers of ${\cal{L}}_K^{1/2}$.

We begin with an estimation for certain integrals involving the generalized Jackson kernels.

\begin{lem}\label{lemmaestimate3} Let $\gamma$ be a positive real number.\ If $2l>\gamma+\alpha(m)+n+1$, then
$$
\int_{0}^{\pi}J_{\nu,n}(t)t^{\gamma}V_m(t) (\sin t)^n\, dt\leq \frac{d_{m,\gamma,l}^{\alpha(m)}}{\mu^{\gamma}},
$$
where $d_{m,\gamma,l}^{\alpha(m)}$ is a constant depending of $m$, $\gamma$, $l$ and the constant $\alpha(m)$ from $A2$.
\end{lem}

\pf The idea of the proof is to detach the normalizing constant $c_{\nu,n}$ from the integral, to find a lower bound for it and an upper bound for the resulting integral.\ Clearly,
$$
c_{\nu,n}\geq\int_{0}^{\pi/2}\left[\frac{\sin(\mu t/ 2)}{\sin(t/2)}\right]^{2l}V_m(t) (\sin t)^n\, dt \geq \frac{2^{2l+n}}{\pi^n}\int_{0}^{\pi/2}t^{n-2l}\left[\sin(\mu t/ 2)\right]^{2l}V_m(t) dt.
$$
From the inequality $V_m(t) \geq c_m t^{\alpha(m)}$, $t \in (0,\pi/2)$, we obtain
$$
c_{\nu,n}\geq c_m\frac{2^{2l+n}}{\pi^n} \int_{0}^{\pi/2}t^{\alpha(m)+n-2l}\left[\sin(\mu t/ 2)\right]^{2l}dt.
$$
The change of variables $s=\mu t$ and the inequality $\mu\geq 2$ provide the estimate
\begin{eqnarray*}
\int_{0}^{\pi/2}t^{\alpha(m)+n-2l}\left[\sin(\mu t/ 2)\right]^{2l}dt&=&\frac{1}{\mu^{\alpha(m)+n+1-2l}}\int_{0}^{\mu\pi/2}s^{\alpha(m)+n-2l}\left[\sin(s/ 2)\right]^{2l}ds
\\ &\geq&\frac{1}{\mu^{\alpha(m)+n+1-2l}}\int_{0}^{\pi}s^{\alpha(m)+n-2l}\left[\sin(s/ 2)\right]^{2l}ds,
\end{eqnarray*}
while an additional adjustment leads to
$$
\int_{0}^{\pi/2}t^{\alpha(m)+n-2l}\left[\sin(\mu t/ 2)\right]^{2l}dt \geq \frac{\pi^{-2l}}{\mu^{\alpha(m)+n+1-2l}}\int_{0}^{\pi}s^{\alpha(m)+n}ds.
$$
The final lower estimate for $c_{\nu,n}$ is
$$
c_{\nu,n} \geq c_m \frac{2^{2l+n}}{\pi^{2l+n}}\frac{1}{\mu^{\alpha(m)+n+1-2l}}\int_{0}^{\pi}s^{\alpha(m)+n}ds=\frac{c_m}{\alpha(m)+n+1} \frac{2^{2l+n}}{\pi^{2l-\alpha(m)-n -1}}\mu^{2l-\alpha(m)-n-1}.
$$
Next, we move to an upper bound for the integral
$$
I:=\int_{0}^{\pi}\left[\frac{\sin(\mu t/ 2)}{\sin(t/2)}\right]^{2l}t^{\gamma}V_m(t) (\sin t)^n\,  dt.
$$
Since $V_m(t)\leq C_mt^{\alpha(m)}$, $t \in (0,\pi)$, it is clear that
$$
I \leq C_m\pi^{2l}\int_{0}^{\pi}[\sin(\mu t/ 2)]^{2l}t^{\gamma+\alpha(m)+n-2l}dt.
$$
Using the change of variables $s=\mu t/2$, we can estimate the integral appearing above as follows
\begin{eqnarray*}
\int_{0}^{\pi}t^{\gamma+\alpha(m)+n-2l}\left[\sin(\mu t/ 2)\right]^{2l}dt &=&\left(\frac{2}{\mu}\right)^{\gamma+\alpha(m)+n+1-2l}\int_{0}^{\mu\pi/2}s^{\gamma+\alpha(m)+n}\left(\frac{\sin s}{s}\right)^{2l}ds\\ &\leq&\left(\frac{2}{\mu}\right)^{\gamma+\alpha(m)+n+1-2l}\int_{0}^{\infty}s^{\gamma+\alpha(m)+n}\left(\frac{\sin s}{s}\right)^{2l} ds.
\end{eqnarray*}
The assumption $2l>\gamma+\alpha(m)+n+1$ guarantees the convergence of the improper integral.\ Proceeding, we have that
$$
I\leq C_m\pi^{2l}\left(\frac{2}{\mu}\right)^{\gamma+\alpha(m)+n+1-2l}\int_{0}^{\infty}s^{\gamma+\alpha(m)+n}\left(\frac{\sin s}{s}\right)^{2l}ds.
$$
Combining our findings, it is promptly seen that the inequality in the statement of the lemma follows and the proof is complete.\eop

\begin{lem}\label{lemmaestimate1} Let $\mathcal{L}_K$ be a positive integral operator generated by a $(T_t,B,\beta)$-H\"{o}lder kernel $K$.\ If $f\in L^2(S^m)$ and $x\in S^m$, then
$$
\left|\mathcal{L}_K^{1/2}(f)(x)-\Phi_{\nu,n}(\mathcal{L}_K^{1/2}(f))(x)\right|\leq \|f\|_2\int_{0}^{\pi}J_{\nu,n}(t)t^{\beta/2}[B(x)+T_t(B)(x)]^{1/2}V_m(t)(\sin t)^n\,dt.
$$
\end{lem}

\pf Fix $f \in L^2(S^m)$ and $x \in S^m$.\ The normalization for the Jackson kernels implies that
$$
\mathcal{L}_K^{1/2}(f)(x)-\Phi_{\nu,n}(\mathcal{L}_K^{1/2}(f))(x) = \int_{0}^{\pi}J_{\nu,n}(t)\left[\mathcal{L}_K^{1/2}(f)(x)-T_t(\mathcal{L}_K^{1/2}(f))(x)\right]V_m(t)(\sin t)^n\, dt.
$$
Hence,
$$
\left|\mathcal{L}_K^{1/2}(f)(x)-\Phi_{\nu,n}(\mathcal{L}_K^{1/2}(f))(x)\right| \leq \int_{0}^{\pi}J_{\nu,n}(t)|D_t(x)|V_m(t)(\sin t)^n\, dt,
$$
where $$D_t(x)=\mathcal{L}_K^{1/2}(f)(x)-T_t(\mathcal{L}_K^{1/2}(f))(x), \quad t\in(0,\pi).$$
The proof will be complete as long as we can reach the estimate below
$$
|D_t(x)|\leq \|f\|_2t^{\beta/2}[B(x)+T_t(B)(x)]^{1/2}, \quad t\in (0,\pi).
$$
Since
$$\mathcal{L}_K^{1/2}(f)(x)=\int_{S^m}K_{1/2}(x,y)f(y)d\sigma_m(y), \quad f \in L^2(S^m),$$
it is easily seen that
$$D_t(x) = \int_{S^m}K_{1/2}(x,y)f(y)d\sigma_m(y) - \frac{1}{\tau_m}\int_{S^m}\int_{S^m}\mathcal{Z}_t^m(x\cdot w)K_{1/2}(w,y)f(y)d\sigma_m(w)d\sigma_m(y), $$
while a change in the integration order leads to
$$D_t(x)=  \frac{1}{\tau_m}\int_{S^m}\left(\tau_m K_{1/2}(x,y) - \int_{S^m}\mathcal{Z}_t^m(x\cdot w)K_{1/2}(w,y)d\sigma_m(w)\right)f(y)d\sigma_m(y).$$
To proceed, we apply H\"{o}lder's inequality to deduce that
$$
|D_t(x)|\leq \frac{1}{\tau_m}\|I_x^t\|_2\|f\|_2, \quad t \in (0,\pi),
$$ in which
$$
I_x^t(y)=\tau_mK_{1/2}(x,y)-\int_{S^m}\mathcal{Z}_t^m(x\cdot w)K_{1/2}(w,y)d\sigma_m(w),\quad y\in S^m,\quad t\in (0,\pi).
$$
The rest of the proof will consist of a tricky estimation of the quantity $\|I_x^t\|_2$ ($t$ fixed).\ A simple calculation leads to
\begin{eqnarray*}
\|I_x^t\|_2^2 & = & \int_{S^m}\tau_m^2K_{1/2}(x,y)K_{1/2}(y,x)d\sigma_m(y) \\
   &  & \hspace*{4mm} + \int_{S^m}\left[-\tau_mK_{1/2}(x,y)\int_{S^m}\mathcal{Z}_t^m(x\cdot w)K_{1/2}(y,w)d\sigma_m(w)\right.\\
&  & \hspace*{6mm} - \tau_mK_{1/2}(y,x)\int_{S^m}\mathcal{Z}_t^m(x\cdot w)K_{1/2}(w,y)d\sigma_m(w)\\
&  & \hspace*{8mm}+ \left.\int_{S^m}\int_{S^m}\mathcal{Z}_t^m(x\cdot w)\mathcal{Z}_t^m(x\cdot z)K_{1/2}(w,y)K_{1/2}(y,z)d\sigma_m(z)d\sigma_m(w)\right]d\sigma_m(y).
\end{eqnarray*}
Interchanging the order of integration and applying the recovery formula (\ref{recover}), we deduce that
\begin{eqnarray*}
\|I_x^t\|_2^2 & = & \tau_m^2K(x,x)-\tau_m\int_{S^m}\mathcal{Z}_t^m(x\cdot w)K(x,w)d\sigma_m(w)\\
&  - &  \tau_m\int_{S^m}\mathcal{Z}_t^m(x\cdot w)K(w,x)d\sigma_m(w)+\int_{S^m}\int_{S^m}\mathcal{Z}_t^m(x\cdot w)\mathcal{Z}_t^m(x \cdot z)K(w,z)d\sigma_m(z)\sigma_m(w).
\end{eqnarray*}
An additional adjustment produces the formula
$$
\|I_x^t\|_2^2\leq  \tau_m^2\left|K(x,x)-T_t(K(x,\cdot))(x)\right| + \tau_m\int_{S^m}\mathcal{Z}_t^m(x\cdot w)\left|K(w,x)-T_t(K(w,\cdot))(x)\right|d\sigma_m(w).
$$
Now, introducing the inequality defining the H\"{o}lder condition, we obtain
$$
\|I_x^t\|_2^2 \leq \tau_m^2 t^{\beta}B(x) + \tau_mt^\beta \int_{S^m}\mathcal{Z}_t^m(x\cdot w)B(w)d\sigma_m(w)=\tau_m^2 t^\beta[B(x) + T_t(B)(x)],
$$
Combining all these findings lead to the inequality in the statement of the lemma.\eop

Next, we not only show that the operator $\mathcal{L}_K^{1/2}-\Phi_{\nu}(\mathcal{L}_K^{1/2}): L^2(S^m)\rightarrow L^{\infty}(S^m)$ is well defined but we also bound the elements on its image.\ Two properties of the norm $\|\cdot \|_{\infty}$ are used in the arguments: the Minkowski's inequality for integrals (\cite[p.194]{folland}) and the inequality $\|\sqrt{f}\|_\infty \leq \|f\|_{\infty}^{1/2}$, $f \in L^\infty(S^m)$, which holds whenever $ f$ is a nonnegative function.

\begin{lem}\label{lemmaestimate2} Let $\mathcal{L}_K$ be a positive integral operator generated by a $(T_t,B,\beta)$-H\"{o}lder kernel $K$.\ If $f\in L^2(S^m)$, then
$$\left\|\mathcal{L}_K^{1/2}(f)-\Phi_{\nu,n}(\mathcal{L}_K^{1/2}(f))\right\|_{\infty}\leq \|B\|_{\infty}^{1/2}(1+M)^{1/2}\left[\int_{0}^{\pi}J_{\nu}(t)t^{\beta/2}V_m(t) (\sin t)^n \,dt\right]\|f\|_2,$$
in which $M$ is a uniform bound for the family $\{\mathcal{Z}_t^m: t\in(0,\pi)\}$ in $L^1([-1,1], d\omega_m)$.
\end{lem}

\pf Fix $f \in L^2(S^m)$ and write
$$
G_{\nu}(f)(x):=\mathcal{L}_K^{1/2}(f)(x)-\Phi_{\nu}(\mathcal{L}_K^{1/2}(f))(x), \quad x \in S^m.
$$
Lemma \ref{lemmaestimate1} and the remarks preceding the lemma imply that
$$
\left|G_{\nu}(f)(x)\right|\leq \|f\|_2\int_{0}^{\pi}J_{\nu,n}(t)t^{\beta/2}\left\|B+T_t(B)\right\|_{\infty}^{1/2}V_m(t) (\sin t)^n \,dt,\quad x \in S^m.
$$
If $M$ is a uniform bound for the family $\{\mathcal{Z}_t^m: t\in(0,\pi)\}$ in $L^1([-1,1], d\omega_m)$, we have that
$$\|T_t(B)\|_{\infty}\leq \|T_t\| \|B\|_{\infty} \leq M\|B\|_{\infty},\quad t \in (0,\pi).$$
The inequality in the statement of the lemma follows.\eop

A similar procedure provides an inequality for the $L^2$-norm of the elements in the image of $G_\nu$.

\begin{lem}\label{corestimatel2} Let $\mathcal{L}_K$ be a positive integral operator generated by a $(T_t,B,\beta)$-H\"{o}lder kernel $K$.\ If $f\in L^2(S^m)$, then
$$\left\|\mathcal{L}_K^{1/2}(f)-\Phi_{\nu,n}(\mathcal{L}_K^{1/2}(f))\right\|_2\leq \|B\|_2^{1/2}(1+M)^{1/2}\left[\int_{0}^{\pi}J_{\nu,n}(t)t^{\beta/2}V_m(t) (\sin t)^n \,dt\right]\|f\|_2,$$
in which $M$ is a uniform bound for the family $\{\mathcal{Z}_t^m: t\in(0,\pi)\}$ in $L^1([-1,1], d\omega_m)$.

\end{lem}


\subsection{Estimates for the approximation numbers of ${\cal{L}}_K$}

In the first step of the subsection, we provide estimates for some of the approximation numbers of ${\cal{L}}_K^{1/2}$ whenever ${\cal{L}}_K$ is generated by a $(T_t,B,\beta)$-H\"{o}lder.\ That implies a decay for the sequence of approximation numbers of ${\cal{L}}_K$ in the case when the corresponding approximating operator has finite rank, details of which are provided in the second step.\ In the last one, we finally obtain decay rates for the sequence of eigenvalues of ${\cal{L}}_K$.
Here, if $l$ is an integer, we will write  $2+l\mathbb{Z}_+:=\{2+l, 2+2l, \ldots\}$.

\begin{thm}\label{thmestimate4} Let $\mathcal{L}_K$ be a positive integral operator generated by a $(T_t,B,\beta)$-H\"{o}lder kernel $K$.\ If there exists a fixed positive integer $q$ so that $\rho(\Phi_{n,k})\leq [q(n+1)]^{\alpha(m)}$, $n=1,2,\ldots$, then
$$a_{(q n)^{\alpha(m)}}(\mathcal{L}_K^{1/2}) =O(n^{-\beta/2}), \quad (n \to \infty).$$
\end{thm}

\pf We begin the proof reminding the reader that the setting and assumptions listed in Section 2 holds here due to the fact that the generating kernel $K$ is $(T_t,B,\beta)$-H\"{o}lder.\ Since the sequence $\{a_{n}(\mathcal{L}_K^{1/2})\}$ is decreasing, if the rank of $\Phi_{n,k}$ has the bound mentioned in the statement of the theorem, it follows that
\begin{eqnarray*}
a_{(qn)^{\alpha(m)}}(\mathcal{L}_K^{1/2})& \leq &a_{(q(n-1))^{\alpha(m)}+1}(\mathcal{L}_K^{1/2})\\
& \leq & \mbox{min}\{\|\mathcal{L}_K^{1/2}-U\|: \rho(U)\leq [q(n-1)]^{\alpha(m)}\}\\
&\leq& \|\mathcal{L}_K^{1/2}-\Phi_{n-2,k}(\mathcal{L}_K^{1/2}):L^2(S^m)\rightarrow L^2(S^m)\|, \quad n=2,3, \ldots.
\end{eqnarray*}
To proceed, choose an integer $l$ in such a way that $2l$ is both, a multiple of ${\alpha(m)}$ and at least $(\beta/2)+\alpha(m)+k+1$.\ If $n\in 2+l\mathbb{Z}_+$, say, $n=2+l\mu$, for some $\mu$, then we can apply Lemma \ref{corestimatel2} (with $\nu=n-2$) to deduce that
$$\|\mathcal{L}_K^{1/2}-\Phi_{n-2,k}(\mathcal{L}_K^{1/2})\| \leq \|B\|_2^{1/2}(1+M)^{1/2}\left(\frac{d_{m,\beta/2,l}^{\alpha(m)}}{\mu^{\beta/2}}\right).$$
Since $l\mu \geq n$, it follows that
$$a_{(qn)^{\alpha(m)}}(\mathcal{L}_K^{1/2}) \leq \|B\|_2^{1/2}(1+M)^{1/2}\left(\frac{d_{m,\beta/2,l}^{\alpha(m)}l^{\beta/2}}{n^{\beta/2}}\right).$$
If $n \in 3+l\mathbb{Z}_+$, the previous inequality implies that
$$a_{(q(n-1))^{\alpha(m)}}(\mathcal{L}_K^{1/2}) \leq \|B\|_2^{1/2}(1+M)^{1/2}\left(\frac{d_{m,\beta/2,l}^{\alpha(m)}l^{\beta/2}}{(n-1)^{\beta/2}}\right).$$
Since $\{a_{n}(\mathcal{L}_K^{1/2})\}$ decreases, a simple calculation leads to
\begin{eqnarray*}
a_{(q n)^{\alpha(m)}}(\mathcal{L}_K^{1/2}) & \leq & \|B\|_2^{1/2}(1+M)^{1/2}\left(1+\frac{1}{n-1}\right)^{\beta/2}\left(\frac{d_{m,\beta/2,l}^{\alpha(m)}l^{\beta/2}}{n^{\beta/2}}\right)\\
& \leq & \|B\|_2^{1/2}(1+M)^{1/2}2^{\beta/2}\left(\frac{d_{m,\beta/2,l}^{\alpha(m)}l^{\beta/2}}{n^{\beta/2}}\right).
\end{eqnarray*}
Inductively, we can infer that
$$
a_{(q n)^{\alpha(m)}}(\mathcal{L}_K^{1/2}) \leq  \|B\|_2^{1/2}(1+M)^{1/2}2^{(k-2)\beta/2}\left(\frac{d_{m,\beta/2,l}^{\alpha(m)}l^{\beta/2}}{n^{\beta/2}}\right), \quad n \in k+l\mathbb{Z}_+,
$$
whenever $k\in \{2,3,\ldots, l+1\}$.\ In other words, except for finitely many positive integers $n$,
$$
a_{(q n)^{\alpha(m)}}(\mathcal{L}_K^{1/2}) \leq  \|B\|_2^{1/2}(1+M)^{1/2}2^{(l-1)\beta/2}\left(\frac{d_{m,\beta/2,l}^{\alpha(m)}l^{\beta/2}}{n^{\beta/2}}\right).
$$
Replacing the constant $\|B\|_2^{1/2}(1+M)^{1/2}2^{(l-1)\beta/2}d_{m,\beta/2,l}^{\alpha(m)}l^{\beta/2}$ with a larger one, if necessary, we reach an inequality in the form
$$
a_{(q n)^{\alpha(m)}}(\mathcal{L}_K^{1/2}) \leq \frac{C}{n^{\beta/2}}, \quad n \in \mathbb{Z}_+,
$$
in which $C$ is a positive constant depending upon $M$, $l$, $B$ and $\beta$.

The main result in this section is this one.

\begin{thm} \label{mainM} Let $\mathcal{L}_K$ be a positive integral operator generated by a $(T_t,B,\beta)$-H\"{o}lder kernel $K$.\ If there exists a fixed positive integer $q$ so that $\rho(\Phi_{n,k})\leq [q(n+1)]^{\alpha(m)}$, $n=1,2,\ldots$, then
$$a_{n}(\mathcal{L}_K) =O(n^{-\beta/\alpha(m)}), \quad (n \to \infty).$$
Further, if the range of $\mathcal{L}_K^{1/2}$ is entirely composed of continuous functions, then the decay can be improved to
$$a_{n}(\mathcal{L}_K) =O(n^{-1-\beta/\alpha(m)}), \quad (n \to \infty).$$
\end{thm}

\pf The first assertion of the theorem is a direct consequence of the previous proposition.\ As for the second one, we need to use a known technique from functional analysis which we now sketch.\ The additional implied by the assumption on the range of $\mathcal{L}_K^{1/2}$ is that the inclusion map $j$ from the space of all continuous functions on $S^m$ to $L^2(S^m)$ is an absolutely 2-summing operator and the composition $j(I-\Phi_{n-2,k})\mathcal{L}_K^{1/2}$ is a Hilbert Schmidt operator, therefore, an absolutely 2-summing operator as well.\ Writing $\pi_2(*)$ to denote the 2-summing norm of a linear operator, we have the following chain of inequalities ($n \in \mathbb{Z}_+$):
\begin{eqnarray*}
(qn)^{\alpha(m)/2}a_{(qn)^{\alpha(m)}}(j(I-\Phi_{n-2,k})\mathcal{L}_K^{1/2}) & \leq & \pi_2(j(I-\Phi_{n-2,k})\mathcal{L}_K^{1/2}) \\
& = & \|j(I-\Phi_{n-2,k})\mathcal{L}_K^{1/2}\|_{HS} \\
& \leq & \pi_2(j)\| \mathcal{L}_K^{1/2}-\Phi_{n-2,k}(\mathcal{L}_K^{1/2}):L^2(S^m)\rightarrow L^{\infty}(S^m)\|.
\end{eqnarray*}
Recalling Lemmas \ref{lemmaestimate2} and \ref{lemmaestimate3}, the previous inequality yields
$$(qn)^{\alpha(m)/2}a_{(qn)^{\alpha(m)}}(j(I-\Phi_{n-2,k})\mathcal{L}_K^{1/2})\leq \pi_2(j) \frac{C}{n^{\beta/2}},$$
in which $C>0$.\ That implies
$$a_{(qn)^{\alpha(m)}}(j(I-\Phi_{n-2,k})\mathcal{L}_K^{1/2}) \leq \frac{C'}{n^{(\beta+\alpha(m))/2}}, \quad n \in \mathbb{Z}_+,$$
for a positive constant $C'$.\ Since the rule that assigns to every operator the sequence of its
approximation numbers is an $s$-scale and $j\mathcal{L}_K^{1/2}=\mathcal{L}_K^{1/2}$, then
the previous inequality implies that
\begin{eqnarray*}
a_{2(qn)^{\alpha(m)}}(\mathcal{L}_K^{1/2}) & \leq & a_{(qn)^{\alpha(m)}}( \mathcal{L}_K^{1/2}-\Phi_{n-2,k}(\mathcal{L}_K^{1/2}))+a_{(qn)^{\alpha(m)} +1}(j\Phi_{n-2,k}\mathcal{L}_K^{1/2})\\
 & = & a_{(qn)^{\alpha(m)}}(\mathcal{L}_K^{1/2}-\Phi_{n-2,k}(\mathcal{L}_K^{1/2})) \\
 & \leq & \frac{C'}{n^{[\beta+\alpha(m)]/2}},\quad n \in \mathbb{Z}_+.
\end{eqnarray*}
The fact that the approximation numbers form a decreasing sequence is all that is needed in order to see that
$$a_{(2qn)^{\alpha(m)}}(\mathcal{L}_K^{1/2}) \leq a_{2(qn)^{\alpha(m)}}(\mathcal{L}_K^{1/2})\leq \frac{C'}{n^{[\beta+\alpha(m)]/2}}, \quad n \in \mathbb{Z}_+.$$
It is now clear that
$$a_n(\mathcal{L}_K^{1/2}) = O(n^{[-\beta-\alpha(m)]/2\alpha(m)}), \quad (n \to \infty),$$
which implies the assertion of the theorem.
\eop

The decay obtained in the abstract setting of the previous theorem has the same structure of that obtained for a positive integral operator generated by a kernel satisfying a standard H\"{o}lder assumption.\ As so,
we conjecture that the decay in the theorem is not improvable.

Let us return now to a positive integral generated by a kernel satisfying an averaged H\"{o}lder condition.\ Recalling Proposition \ref{operatorthm} and the fact that $N(m+1, \nu+1)=O((\nu+1)^m)$, as $n\to \infty$, we can find a positive $q$ so that
$$
\rho(\Phi_{n,1})\leq (q(n+1))^m,\quad n=1,2,\ldots.
$$
This estimate matches the needs in the theorems proved in Section 4.\ Thus, the following result is an immediate consequence.

\begin{thm}  If $\mathcal{L}_K$ is a positive integral operator generated by a kernel $K$ satisfying the averaged H\"{o}lder condition, then
$$
a_{n}(\mathcal{L}_K) =O(n^{-\beta/m}), \quad (n \to \infty).
$$
Further, if the range of $\mathcal{L}_K^{1/2}$ is entirely composed of continuous functions, then the decay can be improved to
$$a_{n}(\mathcal{L}_K) =O(n^{-1-\beta/m}), \quad (n \to \infty).$$
\end{thm}

Due to the remarks at the end of Section 2, a similar theorem holds for a positive integral operator generated by a kernel $K$ satisfying the Stekelov-mean H\"{o}lder condition.\ Details on that will be not included here.

\section{Final remarks}

Most of the concepts and constructions made in this paper can be recovered when we replace the unit sphere with a compact symmetric space of rank 1.\ Indeed, this space is a Riemannian manifold possessing a harmonic analysis structure very similar to that available on the spheres.\ A good source of information on compact symmetric spaces of rank 1, including concepts and results needed in a possible extension of the results proved here, is the survey paper \cite{platonov}.\ We believe the new arguments needed in the detailing of such extension would not justify the writing of an additional paper.

The decay presented in Theorem 4.6 and its corollaries seems to be optimal within the setting considered.\ Restricting ourselves to the two motivational examples of Section 2, we tried for some time to construct a concrete example matching exactly the decay provided by the corresponding results proved in the paper.\ Unfortunately, we were unable to either construct such an example or substantiate optimality.

Recently, we have developed a new technique to deduce sharp decay rates for the sequence of eigenvalues of positive integral operators based on growth and integrability of Fourier coefficients (\cite{jordao1,jordao}).\ This technique allows one to work in an even more general setting, replacing all the arguments involving the usual spherical convolutions with that of spherical convolutions with measures.\ In particular, this approach permits the inclusion of integral operators generated by kernels satisfying H\"{o}lder assumptions defined by families of general multiplier operators.

A final remark concerns the choice $B \in L^\infty(S^m)$ we have made in our definition for the H\"{o}lder assumption.\ On one hand, the restriction is satisfactory because, in relevant concrete cases the function $B$ is, in fact, constant.\ On the other, it may be not.\ However, a more general assumption on $B$, such as $B \in L^1(S^m)$, only appears in purely theoretical results.
%
%

%
%

\vspace*{1.5cm}

\noindent 
Departamento de
Matem\'atica,\\ ICMC-USP - S\~ao Carlos, Caixa Postal 668,\\
13560-970 S\~ao Carlos SP, Brasil\\ E-mails: tjordao@icmc.usp.br; menegatt@icmc.usp.br

\end{document}